\documentclass{article}
\usepackage[utf8]{inputenc}
\usepackage{amsthm}
\usepackage{amsmath, amssymb, amsthm, mathtools}
\usepackage{graphicx}
\usepackage{bbm}
\usepackage{comment}
\usepackage[margin=1.5in]{geometry}
\usepackage{tikz}
\usepackage{mathtools}
\usepackage{authblk}
\usepackage{hyperref}

\DeclareMathOperator{\N}{\mathbf{N}}
\DeclareMathOperator{\U}{\mathcal{U}}
\DeclareMathOperator{\A}{\mathcal{A}}
\DeclareMathOperator{\C}{\mathcal{C}}
\DeclareMathOperator{\D}{\mathcal{D}}

\DeclareMathOperator{\E}{\mathbb{E}}

\DeclareMathOperator{\PP}{\mathbb{P}}
\newcommand{\qle}{\stackrel{?}{<}}

\newtheorem{theorem}{Theorem}[section]

\newtheorem{proposition}[theorem]{Proposition}

\newtheorem{remark}[theorem]{Remark}

\newtheorem{corollary}[theorem]{Corollary}

\theoremstyle{remark}

\theoremstyle{remark}
\newtheorem*{example}{Example}

\newenvironment{customproof}[1][Proof]{%
  \begin{proof}[\textnormal{\textbf{#1}}]%
}{%
  \end{proof}%
}

\title{On a Conjecture on Uniform Group Drawings in the Coupon Collector Problem}
\author[1,2]{Daniel Berend}
\author[2]{Tomer Sher}
\affil[1]{Institute for the Theory of Computing, Ben-Gurion University, Beer Sheva 84105, Israel.}
\affil[2]{Department of Mathematics, Ben-Gurion University, Beer Sheva 84105, Israel.}

\graphicspath{ {./images/} }

\begin{document}

\maketitle

\begin{abstract} 
We address a conjecture of Schilling concerning the optimality of the uniform distribution in the generalized Coupon Collector's Problem (CCP) where, in each round, a subset (package) of $s$ coupons is drawn from a total of $n$ distinct coupons. While the classical CCP (with single-coupon draws) is well understood, the group-draw variant, where packages of size $s$ are drawn, presents new challenges and has applications in areas such as biological network models.

Consider the set of all distributions over the collection of $\binom{n}{s}$ packages of size $s$. Schilling showed that, for $s=n-1$, the uniform distribution yields the minimal expected time for collecting all coupons. She further conjectured that, for $2\le s\le n-2$, the uniform distribution does not yield the minimum. We prove Schilling’s conjecture in full by presenting ``natural'' non-uniform distributions yielding strictly lower expected collection times. Explicit formulas are provided for the expected number of rounds under these and related distributions.\\ 

\end{abstract}

\section{Introduction}


CCP is a problem in combinatorial probability, which can be formulated as follows: There are $n$ distinct coupons, sampled independently with replacement. The coupons arrive one by one, and each has the same probability $\frac{1}{n}$ of being selected at each step. How many drawings, on average, are needed in order to complete a collection of all coupons? We will refer to this version of the problem as the classical version.

The problem can be traced back to de Moivre \cite{tod}. The mathematical model of the problem is compatible with some real-world problems and therefore has potential applications. 

It is well known, and easy to prove, that the expected number of drawings needed to complete the collection is
\begin{align} \label{1}
    n\cdot H_n = n\cdot\left(1+\frac{1}{2}+\cdots+\frac{1}{n-1}+\frac{1}{n}\right),
\end{align}
where $H_n$ is the $n$-th harmonic number.

Coupon subset collection is a generalization of the classical version. Instead of selecting a single coupon each time, a set of $s$ distinct coupons (henceforward a \emph{package}) is selected, where $s$~is an arbitrary fixed integer between 1 and $n-1$. We are interested in the distribution of the number of rounds needed to get all coupons. 
This group-draw variant of the coupon collector problem has also found applications in biological systems. For example, Giannakis et al.~\cite{giannakis2022exchange} model the sharing of mitochondrial DNA in plant cells as a network-based coupon collection process, where encounters between mitochondria correspond to group draws of genetic elements.

Return to the classical case $s=1$. One may inquire how the situation changes if the probabilities of the various coupons are not the same. Boneh and Hofri \cite{boneh} proved that the expected number of draws attains its minimum if and only if all coupons have the same probability $1/n$ of being drawn. Similarly, when drawing packages of size $s\ge 2$, one may consider the case where the various $\binom{n}{s}$ packages have any probabilities. It is natural to inquire whether the uniform distribution again minimizes the expected collection time. The following conjecture of Schilling \cite{schilling} may thus come as a surprise.
\bigskip

\noindent\textbf{Conjecture A.} \emph{For $2\le s \le n-1$, the uniform distribution leads to the minimal expected number of drawings if and only if $s=n-1$.}
\bigskip

She proved the ``if'' part of the conjecture. Note that after the first package is drawn, the probability of drawing the missing coupon is $\frac{1}{n}$. Therefore, the number of rounds after the first is $\mathrm{Geom}\!\left(1 - \frac{1}{n}\right)$-distributed. Chang and Fang \cite{chen} identified the optimal distribution (i.e., the distribution yielding minimal expected collection time) in the case where $n$ is even and $s=n/2$. An optimal distribution in this case is the uniform distribution over any two disjoint packages of $n/2$ coupons each. (The number of drawings after the first round in this case is $\mathrm{Geom}\!\left(\frac{1}{2}\right)$-distributed.) They also proved that the optimal distribution is unique (up to the identity of the coupons in each of the two packages).
They also conjectured what the optimal distribution is for any $s$ dividing $n$
(namely, that it is the decomposition distribution, defined in the next section).

Caron, Hlynka, and McDonald \cite{caron} showed
that when $s=2$, the only case for which the uniform distribution leads to the optimum is $n = 3$, i.e., they proved the ``only if'' part of Schilling's conjecture for $s=2$.

This work is dedicated to proving Schilling's conjecture. Given $n$ and $s$, we will find bounds on the expectation under the uniform distribution and present better distributions (expectation-wise).

In Section \ref{main results} we state our main results. The proofs are presented in Section \ref{proofs}.

\section{The Main Results}\label{main results}
There are several possible approaches towards proving Conjecture A:
\begin{enumerate}
   
    \item Schilling's approach: view the expectation as a real-valued function of $\binom{n}{s}$ (or rather $\binom{n}{s}-1$) variables; the probabilities we assign to the various packages. The conjecture can be proved by showing that $ \Vec{x} = \left(1/\binom{n}{s},\ldots,1/\binom{n}{s}\right)$ is not a global minimum. Schilling suggested, roughly speaking, doing so by showing that the Hessian matrix, evaluated at $\Vec{x}$, is indefinite. That is, the uniform distribution can be neither a local minimum nor a local maximum. She managed to carry out this idea in some special cases. (For a complete description of what she has done, see \cite{schilling}.)
     \item The wishful approach: explicitly define a distribution, prove that it is optimal and that it yields a strictly smaller expected number of drawings than the uniform distribution.
    \item A intermediate approach: present a specific distribution yielding a lower expected number of drawings than does the uniform distribution.
\end{enumerate}
Note that Schilling's approach does not give any hint as to what the optimal distribution is, or even which distributions outperform the uniform distribution. We will stick to the intermediate approach. Specifically, we will construct several explicit distributions that are better than the uniform.  In fact, some appear as strong candidates for actually being optimal. It is interesting to note that the distributions are of different types for $s \ge \left \lfloor{n/2}\right \rfloor$ and $s < \left \lfloor{n/2}\right \rfloor$.

\begin{figure}[!b]
\centering
\includegraphics[width=0.5\textwidth]{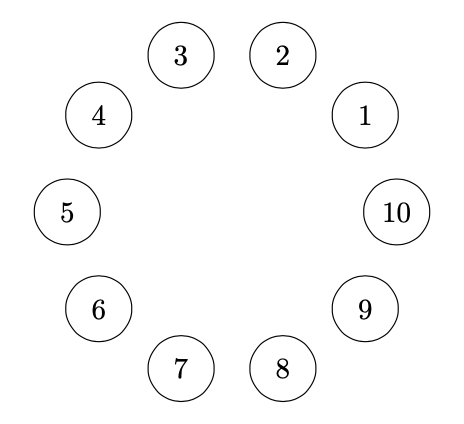}
\caption{Ten coupons along a circle}
\label{arcs}
\end{figure}

The following distributions are important for our work:

\begin{enumerate}
    \item The uniform distribution over all packages of size $s$. We denote this distribution by~$\U_{n,s}$.
    \item The \emph{decomposition distribution}, for $s \mid n$. Here we take the uniform distribution over $n/s$ pairwise disjoint packages of $s$ coupons each. We denote it by $\mathcal{D}_{n,s}$.
    \item The \emph{near-decomposition distribution} for $s \nmid n$. This is an analog of the decomposition distribution.
    Let $q,r \in \N$ be such that $n = q\cdot s + r$ with $0 < r <s$. We divide $n-r$ of the coupons into $q$ pairwise disjoint packages, of size $s$ each. There are $r$ coupons left (the remainder), which we put in a new package. To complete this last package to size $s$, we add to it arbitrary $s-r$ coupons, different from the $r$ already inside. Thus, we have $q+1$ packages in total and take the uniform distribution over them. For simplicity, we denote this distribution also by $\mathcal{D}_{n,s}$, same as the decomposition distribution. In the sequel, when $s \mid n$, we regard the decomposition distribution as a special case of the near-decomposition distribution.
    \item The \emph{arcs distribution}. Given $n$ coupons, enumerated by $1, 2,\ldots, n-1, n$, we take the uniform distribution on the collection of all packages consisting of $s$ ``consecutive'' coupons (modulo $n$), namely
    $$\{1, 2, \ldots, s-1, s\},\ \{2, 3, \ldots, s, s+1\}, \ldots,\ \{n, 1, \ldots, s-2, s-1\}.$$
     It is instructive here to view the coupons as arranged along a circle (see Figure \ref{arcs}). Denote this distribution by~$\mathcal{A}_{n,s}$.

    \item The \emph{complementary decomposition distribution}. We define this distribution only for $s>n/2$. Roughly speaking, we take here the uniform distribution on the complements of all sets taken in the near-decomposition distribution for $n-s$ instead of $s$. More precisely, we do the following, using again arcs. It will be convenient to assume first that $(n-s) \mid n$. We choose the packages as follows. The first package is $\{1,2,\ldots,s\}$. The second is obtained by a rotation of size $n-s$ along the circle, that is $\{n-s+1,n-s+2,\ldots,n\}$. We continue similarly and obtain $m = n/(n-s)$ packages, and take the uniform distribution on this collection. Note that we complete collecting the coupons by drawing any two distinct packages (because $s\ge n/2$). If $(n-s) \nmid n$, we divide into packages similarly, but now we have a ``remainder''. We take $m=\lfloor n/(n-s)\rfloor$ packages, as follows:
    \begin{align*}
        &\{1,2,\ldots,s\}, \{n-s+1,n-s+2,\ldots,n\},\ldots,\\&\{(m-1)\cdot(n-s)+1,(m-1)\cdot(n-s)+2,\ldots,(m-2)(n-s)\}.
    \end{align*}
    (Again, the coupon numbers are to be understood modulo $n$.)
    We denote this distribution by $\C_{n,s}$. 
\end{enumerate}

\begin{example}
For $n=10, s=7$, the sets for the complementary decomposition distribution are $\{1,2,\ldots,7\}, \{4,5,\ldots,10\},\{7,8,\ldots, 3\}.$
\end{example}
For simplicity, we introduce the following
notation for the number of draws required to collect all coupons. This random variable will be denoted by $Y$, with a subscript indicating which of the above distribution types is used. 
Now $Y$ depends on three parameters: 
the distribution type, the total number of coupons, and the package size. 
Accordingly, one should write, for instance, $Y_{\D,10,3}$ to refer to the number of draws needed to obtain all $10$ coupons, drawn in groups of $3$, under the near-decomposition distribution. 
To avoid cumbersome notation, we will simply write $Y_{\U}$, $Y_{\A}$, and so forth, with the understanding that the remaining indices $n$ and $s$ are implicit throughout.

The following theorem presents distributions yielding a lower expectation than the uniform distribution.
\newpage

\begin{theorem} \label{better distributions}
For $n\ge 4$:
\begin{description}
\item{(a)} If $2 \le s < \left \lfloor{n/2}\right \rfloor$, then $\E[Y_{\D}] < \E[Y_{\U}]$.
\item{(b)} If $\left \lfloor{n/2}\right \rfloor \le s \le n-2 $, then $\E[Y_{\A}] < \E[Y_{\U}]$.
\end{description}
\end{theorem}

\begin{corollary}
    Conjecture A is true.
\end{corollary}


The following two propositions provide estimates  or exact expressions for the three quantities appearing in Theorem \ref{better distributions}, namely $\E[Y_{\U}], \E[Y_{\D}]$, and $\E[Y_{\A}]$.

\begin{proposition} \label{expectation bounds}
    For $2 \le s \le n-1$:
    \begin{align*}
        \frac{H_{n}}{H_{n}-H_{n-s}} < \E[Y_{\U}] < 1 + \frac{H_{n-1}}{H_{n}-H_{n-s}}.
    \end{align*}
    In particular, as $n\to\infty$,
     \begin{align*}
         \E[Y_{\U}] =  \frac{H_{n}}{H_{n}-H_{n-s}} + O(1).
     \end{align*}
\end{proposition} 
\noindent

\begin{proposition} \label{better distributions expectations}
For each $n$:
\item{(a)} $\E[Y_{\D}] = \left \lceil{n/s}\right \rceil \cdot H_{\left \lceil{n/s}\right \rceil}$ for $1 \le s \le n$. 
\item{(b)} $\E[Y_{\A}] = \frac{n^2}{s(s+1)}+1$ for $\left \lfloor{n/2}\right \rfloor \le s \le n-1$.
\end{proposition}

\begin{remark} The question regarding the expected time required to obtain all coupons when using the arcs distribution may be considered as a discrete analogue of the following classic problem: Let $C$ be a circle of circumference 1. Draw uniformly randomly arcs of length $a$ on the circle, where $0<a<1$. What is the expected number of drawings required to cover the circle? The second part of the proposition deals with the discrete case, where $s$ corresponds to $a\cdot n$. We are interested here only in $s\ge \lfloor n/2\rfloor$, which corresponds to $a\ge 1/2$, but the continuous problem is interesting and was studied for any $0<a<1$. We refer to \cite{stevens} for the solution of this problem. (See also \cite{solomon1978arcs,WeissteinCircleCovering,PossiblyWrongCouponVariants,CookConsecutiveCoupon}).
\end{remark}

Although Theorem \ref{better distributions} implies Conjecture A, it is interesting to note that we can do even better than the arcs distribution for $s \ge n/2$ if $n$ is sufficiently large using the complementary decomposition distribution. 

Since we complete collecting the coupons by drawing any two distinct packages, after the first drawing, the number of drawings needed to complete a collection is distributed geometrically with a parameter $1-1/m$, where $m=\lfloor n/(n-s)\rfloor$. Therefore,
\begin{equation} \label{formula for w_n}
    \E[Y_{\C}] = 1+\frac{1}{1-1/m} = 1+\frac{1}{1-1/\lfloor n/(n-s)\rfloor}.
\end{equation}
Note that in the case where $(n-s) \mid n$, it readily follows from Proposition \ref{better distributions expectations} that $\E[Y_{\C}] = 1+n/s < \E[Y_{\A}]$ for every $n$. The following proposition refers to the general case where $n$ is not necessarily divisible by $n-s$.

\begin{proposition} \label{even better distribtuion}
For $n > 12$ and $s \ge n/2$, we have $\E[Y_{\C}] < \E[Y_{\A}]$. In particular, for such $n$ and $s$ we have $\E[Y_{\C}] < \E[Y_{\U}]$.
\end{proposition}

\section{Proofs}\label{proofs}
It will be convenient to start with the proof of Proposition \ref{expectation bounds}.
 
\begin{customproof}[Proof of Proposition \ref{expectation bounds}]

The proof proceeds by decomposing the process into mini-steps. Each package draw can be viewed as a sequence of $s$ single-coupon draws. These $s$ drawings have somewhat different statistics. For example, the first coupon may be any coupon, but as we continue, the drawn coupons are known to be distinct from those preceding them within the same step. Consequently, the probability of drawing a previously unseen coupon changes from one mini-step to the next. A precise analysis of this decomposition leads to the desired result.

At each step we collect $s$ distinct coupons. It will be convenient for our purposes to regard each such step as made out of $s$ mini-steps. In each mini-step, a single coupon is drawn such that the coupon that is collected at mini-step $k\cdot s + j$, where $k \ge 0$ and $2 \le j \le s$, is different from all coupons collected at mini-steps $k\cdot s + j'$ for $1 \le j' \le j-1$. However, if we have finished collecting all $n$ coupons after some number $m$ of mini-steps, where $m$ is not divisible by $s$, we must make $s\cdot (1-\{m/s\})$ additional mini-steps. We denote the number of mini-steps needed to obtain all coupons by $Y'$. Clearly, $Y_{\U} =\frac{1}{s} Y'$. 
Denote, for $i \ge 1$,
\begin{align*}
    I_{i} =
    \begin{cases}
        1, &\text{a coupon is still drawn at mini-step $i$}, \\
        0, &\text{otherwise}.
    \end{cases}
\end{align*}
\newline
We have $Y' = \sum_{i=1}^\infty I_{i}$, and in particular $\E[Y'] = \sum_{i=1}^\infty \E[I_{i}]$. By definition, for all $k \ge 0$ we have $I_{k s+1}=I_{ks+2}= \cdots = I_{ks+s}$.

The process described above is very similar to the process in the classical CCP. However, we may finish slightly sooner since in mini-steps of the form $ks+j$ with $2 \le j \le s$, some coupons that already showed up are excluded from being drawn. On the other hand, we may finish somewhat later due to the restriction that the number of mini-steps has to be divisible by $s$.
Now define random variables $(Z_{i})_{i=1}^\infty$, that will let us tie these processes together, by
\begin{align*}
    Z_{i}=
    \begin{cases}
        I_{i}, & \quad i \equiv 1 \pmod{s},\\
        G_{1,i} \cdot I_{i}, & \quad i \equiv 2 \pmod{s},\\
        \vdots & \quad \vdots\\
        G_{s-2,i} \cdot I_{i}, & \quad i \equiv s-1 \pmod{s},\\
        G_{s-1,i} \cdot I_{i}, & \quad i \equiv 0 \pmod{s},
    \end{cases}
\end{align*}

where $G_{j,i}\sim$ Geo$(\frac{n-j}{n})$ is independent of $I_{i}$ for $1 \le j \le s-1$. Let $Z = \sum_{i=1}^\infty Z_{i}$. Then,
\begin{equation} \label{Yexp}
    \begin{alignedat}{2}
    \E[Z] &= \sum_{i=1}^\infty \E[Z_{i}]\\
    &= \sum_{k=0}^\infty \E[I_{ks+1}] + \sum_{k=0}^\infty \E[G_{1,ks +2}\cdot I_{ks +2}] + \cdots + \sum_{k=0}^\infty \E[G_{s-1,ks+s}\cdot I_{ks+s}] \\
    &= \sum_{k=0}^\infty \left(\E[I_{ks+1}]+\frac{n}{n-1}\cdot \E[I_{ks+2}]\right.\\
    &  \quad\quad\quad\,\, \left. +\frac{n}{n-2}\cdot \E[I_{ks+3}]+ \cdots + \frac{n}{n-(s-1)}\cdot \E[I_{ks+s}]\right) \\
    &= \sum_{k=0}^\infty \left(\frac{n}{n}+\frac{n}{n-1}+\frac{n}{n-2}+\cdots+\frac{n}{n-(s-1)} \right)\cdot \E[I_{ks+1}] \\
    &=n\cdot(H_{n}-H_{n-s})\sum_{k=0}^\infty \E[I_{ks+1}] \\ &=n\cdot(H_{n}-H_{n-s})\cdot\frac{1}{s}\sum_{i=1}^\infty \E[I_{i}] \\ &=\frac{n}{s}\cdot(H_{n}-H_{n-s})\E[Y'].
    \end{alignedat}
\end{equation} 

The penultimate equality holds since $I_{k s+1}=I_{ks+2}= \cdots = I_{ks+s}$ for all $k \ge 0$. By the discussion above, $Z$ is very similar to the random variable counting the drawings in the classical CCP. The two would be the same if we were guaranteed that our process terminates after a number of mini-steps divisible by $s$. If the process ends at a mini-step $ks +j$, where $1 \le j \le s-1$, we ``pay'' $s-j$ additional geometric random variables, and therefore more drawings than the classical CCP process. Thus, we obtain a lower and an upper bound for $\E[Z]$:

\begin{align*}
    n\cdot H_{n} < \E[Z] < n&\cdot H_{n} + \sum_{j=1}^{s-1}\frac{n}{n-j} =  n\cdot H_{n} + n\cdot(H_{n-1}-H_{n-s}).
\end{align*}
By ($\ref{Yexp}$):
\begin{align*}
     \frac{s\cdot H_{n}}{H_{n}-H_{n-s}} < \E[Y'] < \frac{s\cdot(H_{n} -H_{n-s}+H_{n-1})}{H_{n}-H_{n-s}}.
\end{align*}
Since $Y_{\U} =\frac{1}{s}Y'$,
\begin{align*}
     \frac{H_{n}}{H_{n}-H_{n-s}} < \E[Y_{\U}] < \frac{H_{n} -H_{n-s}+H_{n-1}}{H_{n}-H_{n-s}} =  1 + \frac{H_{n-1}}{H_{n}-H_{n-s}}.
\end{align*}
\end{customproof}

\begin{customproof}[Proof of Proposition \ref{better distributions expectations}]
(a) In the near-decomposition distribution, the coupons are split into $\left \lceil{n/s}\right \rceil$ packages with equal probability. We complete a collection when all $\left \lceil{n/s}\right \rceil$ packages have been obtained. That is, it reduces to the classical CCP with $\left \lceil{n/s}\right \rceil$ coupons, and by (\ref{1}) the expectation is
$$\left \lceil{n/s}\right \rceil \cdot H_{\left \lceil{n/s}\right \rceil}.$$
(b) For $0 \le k \le n-s$ and $k=n$, denote by $t_{k}$ the expected number of drawings left to collect $k$ missing coupons after we have already collected $n-k$. Clearly, $t_n = \E[Y_{\A}]$. Note that $t_{k}$ is well-defined because, at any step, all $k$ missing coupons are consecutive, i.e., there is a missing arc of $k$ coupons. This holds by symmetry and the fact that $s \ge \left \lfloor{n/2}\right \rfloor$. Suppose that at some point there are still $k$ missing coupons. Denote by $A_{i}$ the event whereby, at the next drawing of $s$ coupons, we obtain $i$ of the missing, $0 \le i \le k$. By the law of total expectation,
\begin{align*}
   t_{k}=1+\PP(A_{0})t_{k}+\PP(A_{1})t_{k-1}+\cdots+\PP(A_{k-1})t_{1}+\PP(A_{k})t_{0}, \quad 1 \le k \le n-s,
\end{align*}
with the initial condition $t_{0}=0$. This yields
\begin{equation}\label{recursion of t}
   t_{k} = \frac{1}{1-\PP(A_{0})}\cdot\left(1+\PP(A_{1})t_{k-1}+\cdots +\PP(A_{k-1})t_{1}+\PP(A_{k})t_{0}\right).
\end{equation}
Continuing to view the $n$ coupons as arranged along a circle (see Fig. \ref{arcs}), we conclude that, since $1 \le k \le n-s$
\begin{equation} \label{P_a}
    \PP(A_{i}) =
    \begin{cases}
        (n-(s+k-1))/n, &\text{$i=0$}, \\
        2/n, &\text{$1 \le i \le k-1$},\\
        (s-k+1)/n, &\text{$i=k$}.
    \end{cases}
\end{equation}
By (\ref{P_a}), we may rewrite (\ref{recursion of t}) in the form
\begin{equation} \label{formula of t_k}
    t_{k}= \frac{n}{s+k-1}\left( 1+\sum_{i=1}^{k-1}\frac{2}{n}\cdot t_{k-i}\right)=\frac{n}{s+k-1}+\frac{2}{s+k-1}\sum_{i=1}^{k-1} t_{i}.
\end{equation}
We now prove by induction that
\begin{equation}\label{9}
    t_{k} = \frac{n(s+k)}{s(s+1)} ,\quad\text{ $1 \le k \le n-s$.}
\end{equation}

For $k=1$, this is immediate from (\ref{formula of t_k}). Now assume that (\ref{9}) holds up to some $k$, and prove it for $k+1$. By (\ref{formula of t_k})
\begin{align*}
    t_{k+1} &=\frac{n}{s+k}+\frac{2}{s+k}\sum_{i=1}^{k} t_{i} = \frac{n}{s+k}+\frac{2}{s+k}\sum_{i=1}^{k} \frac{n(s+i)}{s(s+1)} \\&=\frac{n}{s+k}+\frac{2n}{s(s+1)(s+k)}\sum_{i=1}^{k} (s+i) = \frac{n}{s+k}+\frac{2n}{s(s+1)(s+k)}\cdot\frac{k(s+1+s+k)}{2} \\&= \frac{n(s+k+1)}{s(s+1)}. 
\end{align*}
In particular, (\ref{9}) yields $$t_{n-s} = \frac{n^2}{s(s+1)}.$$ Since $ t_{n} = 1 +t_{n-s}$, we get
\begin{equation}\label{formula for t_n}
    \E[Y_{\A}] = 1+\frac{n^2}{s(s+1)} ,\quad\text{$s \ge \left \lfloor{n/2}\right \rfloor$.}
\end{equation}
\end{customproof}

\begin{customproof}[Proof of Theorem \ref{better distributions}]

The proof compares the expected collection times under the relevant distributions by analyzing inequalities involving harmonic sums. For small $s$, the near-decomposition distribution yields smaller expectations than the uniform one, shown via analytic bounds on harmonic numbers; for larger $s$, a recursive argument on the expected number of missing coupons proves that the arcs distribution is superior.

(a) By Proposition \ref{expectation bounds}
\begin{align*}
    \E[Y_{\U}] > \frac{H_{n}}{H_{n}-H_{n-s}},
\end{align*}
so that, by Proposition \ref{better distributions expectations}.(a), it suffices to show that
\begin{align} \label{6}
    \left \lceil{n/s}\right \rceil \cdot H_{\left \lceil{n/s}\right \rceil}  < \frac{H_{n}}{H_{n}-H_{n-s}}.
\end{align}
We distinguish between the case where $s \mid n$ and the case where $s \nmid n$. First we assume $s \mid n$. We have:
\begin{align*}
    \frac{H_{n}}{H_{n}-H_{n-s}} &= \frac{\sum_{j=1}^{n/s} \sum_{i=1}^s \frac{1}{(j-1)s+i}}{\sum_{i=1}^s \frac{1}{n-s+i}} = \sum_{j=1}^{n/s} \left(\frac{\sum_{i=1}^s \frac{1}{(j-1)s+i}}{\sum_{i=1}^s \frac{1}{n-s+i}}\right).
\end{align*}

In general, for arbitrary positive numbers $a_i$ and $b_i$, $1 \le i \le m$, we have $$\frac{\sum_{i=1}^m a_i}{\sum_{i=1}^m b_i} = \sum_{i=1}^m \left(\frac{b_i}{\sum_{j=1}^m b_j}\right) \frac{a_i}{b_i} \ge \min_{1 \le i \le m} \frac{a_i}{b_i}$$
and the inequality is sharp unless all ratios $a_i / b_i$ are equal.

Hence for each $j$:

\begin{align*}
    \frac{\sum_{i=1}^s \frac{1}{(j-1)s+i}}{\sum_{i=1}^s \frac{1}{n-s+i}} > \min_{1 \le i \le s}\frac{n-s+i}{(j-1)s+i} = \min_{1 \le i \le s}\frac{n-js}{(j-1)s+i} + 1 = \frac{n}{js}.
\end{align*}
Therefore:

\begin{align*}
    \frac{H_{n}}{H_{n}-H_{n-s}} &> 
     \sum_{j=1}^{n/s} \frac{n}{js} =\frac{n}{s}H_{n/s}.
\end{align*}

Now let $s \nmid n$. There are $q,r \in \N$ such that $n = q\cdot s + r$ and $1 \le r \le s-1$. Since $s <  \left \lfloor{n/2}\right \rfloor $, we have $q \ge 2$. In (\ref{6}), notice that, for a fixed $q$, the left-hand side is constant. The right-hand side is an increasing function of $n$, so we prove (\ref{6}) for $n$ of the form $n = qr+1$. We will use the well-known estimates 
$$ H_{n}-1 \le \log n \le H_{n-1}, $$
which imply 
\begin{align} \label{h_n bounds}
    \log n + \frac{1}{n} \le H_{n} \le \log n + 1.
\end{align}
Note that $q+1 =   \left \lceil{n/s}\right \rceil = \frac{n+s-1}{s}$ and $H_{n}-H_{n-s}=1/n+\cdots+1/(n-s+1) < s/(n-s+1)$, so that it suffices to show that
\begin{align} \label{suffices to show}
    \frac{n+s-1}{s}\cdot H_{q+1}\cdot\frac{s}{n-s+1} < H_{n}.
\end{align}
Indeed,
\begin{align*}
    \frac{n+s-1}{s}\cdot H_{q+1}\cdot\frac{s}{n-s+1} < H_{n} &\iff \frac{qs+1+s-1}{qs+1-s+1}\cdot H_{q+1} < H_{qs+1} \\&\iff \frac{(q+1)s}{(q-1)s+2}\cdot H_{q+1} < H_{qs+1} \\&\iff \frac{q+1}{q-1+2/s}\cdot H_{q+1} < H_{qs+1} \\&\iff \frac{1}{q-1+2/s}+\frac{q+1}{q-1+2/s}\cdot H_{q} < H_{qs+1}.
\end{align*}
We want to find out when the latter inequality holds. Using (\ref{h_n bounds}), and the fact that $q \ge 2$, we can show routinely that it is enough to check when
\begin{align} \label{7}
    1+ \frac{3+2\log q}{q-1} < \log s.
\end{align}
Thus, the correctness of (\ref{7})
will imply the correctness of this part. Some calculation yields that, for $s \ge 220$, inequality (\ref{7}) holds (independently of $q$). For $3 \le s<220$, inequality (\ref{7}) is a ``trade-off inequality'' between $s$ and $q$.
That is, for any value of $s$ we choose, we get a value of $q$ such that our claim holds for every $q' \ge q$ (because the left-hand side of (\ref{7}) is decreasing for $q\ge 2$). We conclude that for $3 \le s<220$ we have a finite number of pairs $(n,s)$ that do not satisfy (\ref{7}), but can easily be checked using computer algebra. Note that, for $s=2$, inequality (\ref{7}) does not hold for any $q$ (as the left-hand side is greater than 1 and the right-hand side is less than 1), so we deal with this case separately.

If $s=2$ and $s \nmid n$ as assumed before, let $k =(n-1)/2$. We are interested in the case where $k \ge 3$. By (\ref{suffices to show}) it is sufficient to prove that
\begin{align*}
    \frac{1}{k}\cdot H_{k+1} < H_{2k+1} - H_{k+1}.
\end{align*}
This inequality holds for $k=5$ and it can easily be checked that the left-hand side is decreasing while the right-hand side is increasing. For $k=3,4$ we verify (\ref{6}) directly.

(b) We first develop a recursion formula for the expected number of drawings under the uniform distribution. Let $a_{k}$ be the expected number of drawings required to collect $k$ missing coupons after we have already collected $n-k$. Note that $a_n = \E[Y_{\U}]$. Let $B_{i}$, $0 \le i \le k$, be the event whereby $i$ new coupons are collected in a single draw while there are still $k$ missing coupons. Clearly
\begin{equation}\label{prob for a}
    \PP(B_{i}) = \frac{\binom{k}{i}\binom{n-k}{s-i}}{\binom{n}{s}}.
\end{equation} Similarly to what we have done with $t_{k}$ in the proof of Proposition \ref{better distributions expectations}(b), we obtain with the initial condition $a_{0}=0$,
\begin{equation} \label{recursion of a}
   a_{k} = \frac{1}{1-\PP(B_{0})}\cdot\left(1+\PP(B_{1})a_{k-1}+\cdots +\PP(B_{s-1})a_{k-(s-1)}+\PP(B_{s})a_{k-s}\right).
\end{equation}

We will prove by induction that $t_{k} < a_{k}$ for $2 \le k \le n-s$.\\
For $k=2$ we have, by (\ref{9}), (\ref{prob for a}), (\ref{recursion of a}) and some routine calculations, $$t_{2} = \frac{n(s+2)}{s(s+1)},\qquad a_{2} = \frac{n(3n-2s-1)}{s(2n-s-1)},$$ so that
\begin{align*}
    t_{2} < a_{2} &\iff \frac{n}{s}\cdot\frac{s+2}{s+1} < \frac{n}{s}\cdot\frac{3n-2s-1}{2n-s-1} \iff 1+\frac{1}{s+1} < 1+\frac{n-s}{2n-s-1}\\ &\iff 2n-s-1 < (n-s)(s+1) \iff (s-1)(s-(n-1)) < 0,
\end{align*}
which clearly holds. Now assume that $t_{i} < a_{i}$ for every $2 \le i \le k-1$. Let us prove that $t_{k} < a_{k}$. By (\ref{recursion of a})
\begin{align*}
    a_{k} &= \frac{1}{1-\PP(B_{0})}\left(1+\sum_{j=1}^{k-1}\PP(B_{j})a_{k-j}\right) > \frac{1}{1-\PP(B_{0})}\left(1+\sum_{j=1}^{k-1}\PP(B_{j})t_{k-j}\right)\\ &= \frac{1}{1-\PP(B_{0})}\left(1+\frac{n}{s(s+1)}\sum_{j=1}^{k-1}(s+k-j)\PP(B_{j})\right) \\ &= \frac{1}{1-\PP(B_{0})}\cdot\frac{n}{s(s+1)}\left(\frac{s(s+1)}{n}+\sum_{j=0}^{k}[(s+k-j)\PP(B_{j})]-(s+k)\PP(B_{0})-s\PP(B_{k})\right) \\ &= \frac{1}{1-\PP(B_{0})}\cdot\frac{n}{s(s+1)}\left(\frac{s(s+1)}{n}+(s+k)\sum_{j=0}^{k}\PP(B_{j})-\sum_{j=0}^k j\PP(B_{j}) -(s+k)\PP(B_{0})-s\PP(B_{k})\right) \\ &= \frac{1}{1-\PP(B_{0})}\cdot\frac{n}{s(s+1)}\left(\frac{s(s+1)}{n}+(s+k)\cdot(1-\PP(B_{0}))-\frac{ks}{n}-s\PP(B_{k})\right) \\ &= \frac{n(s+k)}{s(s+1)} + \frac{1}{1-\PP(B_{0})}\cdot\frac{n}{s+1}\left(\frac{(s+1-k)}{n}-\PP(B_{k})\right) \\ &= t_{k} + \frac{n}{(1-\PP(B_{0}))(s+1)}\cdot(\PP(A_{k})-\PP(B_{k})) \ge t_{k}
\end{align*}
In particular, $t_{n} = 1+t_{n-s} <  1+a_{n-s}= a_{n}$. 
\end{customproof}

\begin{customproof}[Proof of Proposition \ref{even better distribtuion}]
Let $m=\lfloor n/(n-s) \rfloor$. We easily see that: $$\frac{m-1}{m}\cdot n \leq s < \frac{m}{m+1}\cdot n.$$ 
Therefore, if we fix $n$, $\E[Y_{\C,n,s}]$ is constant as a function of $s$ on intervals of the form \newline $[\frac{m-1}{m}\cdot n, \frac{m}{m+1}\cdot n)$. Note that $m \ge 2$ since $s\ge n/2$. Since $t_{n}(s)$ is decreasing as a function of $s$, it is enough to prove that
\begin{align} \label{inequality for 3.2}
    w_n\left(\frac{m-1}{m}\cdot n\right) < t_{n}\left(\frac{m}{m+1}\cdot n\right),\qquad m \ge 2.
\end{align}
(Here, we view $t_n(\cdot)$ and $w_n(\cdot)$ as defined by (\ref{formula for t_n}) and (\ref{formula for w_n}), respectively, for any real argument.) Since $s \le n-2$, we have
\begin{align*}
     m \le \frac{n}{n-s} \le \frac{n}{2}.
\end{align*}
Now we check when (\ref{inequality for 3.2}) holds for $2 \le m \le \frac{n}{2}$:
 \begin{align*}
     w_n\left(\frac{m-1}{m}\cdot n\right) < t_{n}\left(\frac{m}{m+1}\cdot n\right) &\iff 1+\frac{m}{m-1} < 1+n^2\cdot\frac{(m+1)^2}{nm(nm+m+1)} \\&\iff \frac{m}{m-1} < \frac{n(m+1)^2}{nm^2+m^2+m}
     \\&\iff m^3\cdot n +m^3+m^2 < n(m+1)^2(m-1) \\&\iff n > \frac{m^2(m+1)}{m^2-m-1}. 
\end{align*}
Finding the maximum of the right-hand side over $[2,\frac{n}{2}]$, we see that this inequality holds for $n>12$. 
The ``In particular'' part follows immediately from Theorem \ref{better distributions}.
\end{customproof}

\section{Acknowledgments}
We thank the anonymous reviewers for their valuable comments and suggestions that helped improve this work.

\clearpage

\bibliographystyle{plain}
\bibliography{references_main}
\end{document}